\documentclass{osa-article}

\journal{osajournal}

\articletype{Research Article}

\newcommand{\bx}{\mathbf{x}}
\newcommand{\ee}{\mbox{e}}

\newcommand{\be}{\mathbf{e}}   
\newcommand{\bj}{\mathbf{j}}      
\newcommand{\br}{\mathbf{r}}      
\newcommand{\sd}{\text{d}}          
\usepackage{multirow,comment}
\usepackage{pgf,tikz,tikz-3dplot}
\usepackage{pdflscape}
\usepackage{pgfplots,xcolor} 
\usepackage{subfigure,mathtools}
\newcommand{\mc}{\multicolumn{1}{c}}
\DeclarePairedDelimiter{\ceil}{\lceil}{\rceil}

\begin{document}

\title{Circulant preconditioning in the volume integral equation method for silicon photonics}

\author{Samuel P. Groth,\authormark{1,2,*} Athanasios G. Polimeridis,\authormark{3} Alexandra Tambova,\authormark{4} and Jacob K. White\authormark{2}}

\address{\authormark{1}Department of Engineering, University of Cambridge, Cambridge, UK}
\address{\authormark{2}Department of Electrical Engineering and Computer Science, Massachusetts Institute of Technology, Cambridge, MA, USA\\
\authormark{3}Q Bio, Redwood, CA 94063, USA\\
\authormark{4}Skolkovo Institute of Technology, Moscow 143026, Russia}

\email{\authormark{*}samuelpgroth@gmail.com} 


\begin{abstract}
Recently, the volume integral equation (VIE) approach has been proposed as an efficient simulation tool for silicon photonics applications [J.~Lightw.~Technol.~{\bfseries 36}, 3765 (2018)]. However, for the high-frequency and strong contrast problems arising in photonics, the convergence of iterative solvers for the solution of the linear system can be extremely slow. The uniform discretization of the volume integral operator leads to a three-level Toeplitz matrix, which is well suited to preconditioning via its circulant approximation. In this paper, we describe an effective circulant preconditioning strategy based on the multi-level circulant preconditioner of Chan and Olkin [Numer.~Algorithms {\bfseries 6}, 89 (1994)]. We show that this approach proves ideal in the canonical photonics problem of propagation within a uniform waveguide, in which the flow is unidirectional. For more complex photonics structures, such as Bragg gratings, directional couplers, and disk resonators, we generalize our preconditioning strategy via geometrical partitioning (leading to a block-diagonal circulant preconditioner) and homogenization (for inhomogeneous structures). Finally, we introduce a novel memory reduction technique enabling the preconditioner's memory footprint to remain manageable, even for extremely long structures. The range of numerical results we present demonstrates that the preconditioned VIE is fast and has great utility for the numerical exploration of prototype photonics devices. 
\end{abstract}


\section{Introduction}
\label{sec:intro}
Volume integral equation (VIE) methods are popular for wave propagation and scattering problems in which the scattering body is potentially inhomogeneous. Example scenarios include birefringent mineral dust and ice crystals \cite{nousiainen2009single,yurkin2007discrete}, human body models \cite{polimeridis2014stable}, and, more recently, infinite silicon photonics structures truncated by absorbing regions \cite{tambova2018adiabatic}. Discretizing the VIE on a uniform (``voxelized'') grid leads to a three-level Toeplitz system matrix with which a matrix-vector product can be calculated in $\mathcal{O}(n\log n)$ operations via the fast-Fourier transform (FFT), where $n$ is the number of voxels in the grid. Therefore the cost of solving these equations via an iterative Krylov subspace technique, such as the generalized minimum residual method (GMRES), is $\mathcal{O}(mn\log n)$, where $m$ is the number of iterations required for convergence. 

For scenarios in which the body is small relative to the wavelength, the factor $m$ is typically small, leading to very efficient simulations. However, as we shall see later, $m$ can increase dramatically with the increasing relative scale of the problem. In many of the aforementioned problems, the length of the structures can be hundreds or thousands of wavelengths, for which $m$ is so large as to preclude practical VIE simulations. In order to tame this growth in iteration count, an effective preconditioning strategy is required. 

This preconditioning requirement is not unique to VIE, but rather is encountered in all numerical methods for high-frequency wave problems (save for niche high-frequency methods, e.g.~\cite{groth2018hybrid}, which are not appropriate in the present setting). Previous preconditioning strategies for integral equations for different wave scattering scenarios include multiplicative Calder\'{o}n preconditioners~\cite{christiansen2002preconditioner,bagci2010calderon,kleanthous2019calderon} and the inverse fast-multipole method \cite{coulier2017inverse}. However, no preconditioning approach has yet to be studied for integral equations in the extremely high-frequency context of three-dimensional silicon photonics problems. 

In our voxel-based VIE setting, it seems natural to exploit the literature on circulant-based preconditioners for Toeplitz systems, which goes back to the work of T.\ Chan \cite{chan1988optimal}, Strang \cite{strang1986proposal}, Tyrtyshnikov~\cite{tyrtyshnikov1992optimal}, and others in the 1980s and 1990s. The vast majority of the applications employing such preconditioners are in image reconstruction, signal processing, least squares problems, and ordinary differential equations \cite{ng2004iterative}, where the system matrix is purely Toeplitz or Toeplitz-block. The application of circulant preconditioners to block-Toeplitz Toeplitz-block structures, such as the present one, has proved less fruitful since it was shown that using repeated circulant approximations on multiple levels yields increasingly poor preconditioners \cite{capizzano2000any}. 

In this article, we focus on the particular high-frequency problems occurring in silicon photonics applications, in which the structures of interest are typically combinations and deformations of long, thin and shallow waveguides. As such, the high-frequency nature of the problems is confined to one, or at most two, of the three physical dimensions. This suggests that a circulant preconditioning strategy may be effective since we only need to employ the circulant approximation in the one, or two, dimensions of extreme length. We show that, when employed in one dimension for homogeneous structures, this strategy yields small Krylov iteration counts, independent of the  length of the structure, while in two dimensions, the iteration count still grows with length/width but is nevertheless greatly reduced. Further, we present simple strategies to generalize this approach to inhomogeneous structures such as Bragg gratings or when adiabatic absorbers are appended to the domain truncation sites~\cite{tambova2018adiabatic}. Finally, we consider silicon photonics components that are composed of multiple separate structures, such as disk resonators and directional couplers. In such cases, we propose an effective blocked-circulant preconditioner based on partitioning the geometry into boxes, each enclosing a constituent structure. Each block in the blocked-circulant preconditioner is the circulant preconditioner for the corresponding box in the partitioned geometry. 

We remark that a similar preconditioning strategy was considered by Remis \cite{remis2012circulant,remis2013preconditioning} for VIE for acoustic scattering by inhomogeneous cuboids in one and two dimensions. There a generalization of Chan's 1-level optimal preconditioner \cite{chan1988optimal} to inhomogeneous structures was presented and was shown to yield rapid convergence rates for iterative solvers for the simple problems considered. The present work represents the first application of multi-level circulant preconditioners to three-dimension EM problems of practical interest. Furthermore, we present a more pragmatic and faster approach than that of Remis for inhomogeneous structures. Where Remis constructed the optimal circulant approximation to the system matrix, we first average the permittivities of the inhomogeneous structure and construct the traditional Chan-circulant approximation to the resulting system matrix and employ this. For the mostly homogeneous photonics structures considered here, we observe that this approach is fast and effective. A very brief outline of some of our results has been previously presented by the current authors in \cite{Groth2018circulant}.

The layout of the paper is as follows. Section~\ref{sec:VIE} gives a brief review of the VIE formulation and its voxel-based discretization, leading to a block-Toeplitz Toeplitz-block (BTTB) matrix system. In Section~\ref{sec:straight} we examine the performance of GMRES for unpreconditioned-VIE for the canonical photonics problem of wave propagation within a straight waveguide. Here we observe that the iteration count of GMRES increases linearly with waveguide length and also increases with permittivity, showing that VIE for high-frequency high-contrast problems of silicon photonics requires an effective preconditioner in order to be a viable solver. Section~\ref{sec:prec} provides a brief summary of the 1- and 2-level circulant preconditioners employed in this article, originally proposed by T. Chan and Olkin in \cite{chan1994circulant}. In Section~\ref{sec:results} we present numerical results to demonstrate the performance of the circulant preconditioner. In particular, we examine the performance for a straight waveguide, a Bragg grating, a disk resonator, and a directional coupler. In Section~\ref{sec:conc} we offer some conclusions and discuss potential future improvements and generalizations of the method. 

\section{Volume integral equations}
\label{sec:VIE}
Consider the scattering of time-harmonic electromagnetic waves with angular frequency $\omega$ by a non-magnetic, dielectric, potentially inhomogeneous, object occupying a bounded domain $\Omega$ in 3D space $\mathbb{R}^3$. Throughout the time-dependence $\ee^{j\omega t}$ is assumed with $j=\sqrt{-1}$. The electric properties are defined as
\begin{equation}
	 \epsilon = \epsilon_0,\ \mu = \mu_0 \quad \text{in}\ \mathbb{R}^3\backslash\Omega; \quad \epsilon = \epsilon_{r}(\br)\epsilon_0,\ \mu=\mu_0 \quad \text{in}\ \Omega, 
\end{equation}
where $\epsilon_0$ and $\mu_0$ are the free-space permittivity and permeability, respectively. The relative permittivity is written
\begin{equation}
	\epsilon_r(\br)  = \epsilon_r'(\br) - j\epsilon_r''(\br); \quad \epsilon_r(\br)'\in(0,\infty), \epsilon_r''(\br)\in[0,\infty].
\end{equation}
The total electric field $\be$ can be expressed in terms of polarization currents $\bj$ as follows~\cite{chew2001fast}:
\begin{equation}
	\be = \be_{\text{inc}} + \frac{1}{j\omega\epsilon_0}(\mathcal{N}-\mathcal{I})\bj,
	\label{eqn:field_decomp}
\end{equation}
where $\be_{\text{inc}}$ is the incident field (in this paper, a dipole source), $\mathcal{I}$ is the identity operator, and the integro-differential operator $\mathcal{N}$ is defined as
\begin{align}
	\mathcal{N}f &:= \nabla\times\nabla\times\int_{\Omega} \frac{\ee^{-j k_0|\br-\br'|}}{4\pi |\br-\br'|}f(\br')\sd \br' ,
\end{align}
where $k_0=\omega\sqrt{\epsilon_0\mu_0}$ the free-space wavenumber. Also, the equivalent current density is defined in terms of the electric field as
\begin{equation}
	\bj(\br) = j\omega\epsilon_0(\epsilon_r(\br)-1)\be(\br).
	\label{eqn:equiv_currents}
\end{equation}
Combining (\ref{eqn:field_decomp}) and (\ref{eqn:equiv_currents}), one can derive the current-based VIE
\begin{equation}
	(\mathcal{I}-\mathcal{M}\mathcal{N})\bj = j\omega\epsilon_0\mathcal{M}\be_{\text{inc}},
	\label{eqn:JVIE}
\end{equation}
where $\mathcal{M}(\br):=(\epsilon_r(\br)-1)/\epsilon_r(\br)$ is a local medium operator.
Following \cite{polimeridis2014stable}, we term the particular formulation (\ref{eqn:JVIE}) the JVIE.

There are numerous discretization techniques available for numerically solving the JVIE (\ref{eqn:JVIE}). Here we employ the Galerkin method over a uniform (``voxelized'') discretization of the domain. We begin the discretization by choosing an appropriate voxel dimension $\Delta$. In this article we take $\Delta\approx 22/\lambda_{\text{int}}$ to ensure accurate simulations (in line with our previous article~\cite{tambova2018adiabatic}), where $\lambda_{\text{int}}$ is the interior wavelength.
Then a box bounding the scatterer is constructed, of dimension $n_x\Delta \times n_y\Delta \times n_z\Delta$ so that the voxel grid consists of $N = n_x\times n_y\times n_z$ voxels.

We then approximate the unknown currents over the voxel grid as 
\begin{equation}
    \mathbf{j}(\br) \approx \sum_{i=1}^lw^x_i\mathbf{p}_i^x + \sum_{i=1}^mw^y_i\mathbf{p}_i^y + \sum_{i=1}^nw^z_i\mathbf{p}_i^z,
\end{equation}
where, for $\alpha=x,y,z$, the weights $w_i^{\alpha}$ are to be determined and
\[
\mathbf{p}_i^{\alpha}(\br)=
    \begin{cases}
        \frac{1}{\sqrt{V}} \hat{\boldsymbol{\alpha}},& \br\ \text{in voxel $i$,} \\
        0, &\text{otherwise},
    \end{cases}
\]
with $\hat{\boldsymbol{\alpha}} =\hat{\mathbf{x}},\hat{\mathbf{y}},\hat{\mathbf{z}}$ being the unit direction vectors in the $x$-, $y$-, $z$-directions, respectively.
The scaling by the square root of the voxel volume $V$ is included so that $\langle \mathbf{p}^{\alpha}_i,\mathbf{p}^{\beta}_j \rangle = \delta_{\alpha\beta ij}$, where $\langle\cdot,\cdot\rangle$ is the standard $L^2$ inner product and $\delta_{\alpha\beta ij}$ is the generalized Kronecker delta. We also assume that the permittivity $\epsilon_r(\br)$ is piecewise constant, with its value being defined at the center of each voxel.

Applying the Galerkin method to the JVIE (\ref{eqn:JVIE}), with testing functions $\mathbf{p}_i$, gives rise to the linear system of 3N equations:
\begin{equation}
	(\textbf{I}-\textbf{M}\textbf{N})\textbf{w} = j\omega\epsilon_0\textbf{M}\textbf{e}_{\text{inc}},
	\label{eqn:discrete_form}
\end{equation}
where $\textbf{I}$ is the identity matrix, $\textbf{M}$ is a diagonal matrix with the entries corresponding $\mathcal{M}(\mathbf{r})$ evaluated at voxel centers, $\textbf{e}_{\text{inc}}$ is a vector with entries corresponding to the incident field evaluated at voxels centers, and 
\begin{equation}
	\textbf{N}^{\alpha\beta}_{ij} = \langle \mathcal{N}\mathbf{p}^{\beta}_j,\mathbf{p}^{\alpha}_i\rangle.
    \label{eqn:op_N}
\end{equation}
Note that we have assumed that the both the material properties and incident field piecewise constant functions on the voxel grid.
The discrete system (\ref{eqn:discrete_form}) has the form
\begin{equation}
\left[
\textbf{I}-
  \renewcommand{\arraystretch}{0.9}
  \left(
  \begin{array}{ c c c | c c c | c c c }
     & &  &  &  & \mc{}& & &  \\
     & &  &  &  & \mc{} & & &  \\
    \multicolumn{3}{c|}{\raisebox{1\normalbaselineskip}[0pt][0pt]{\footnotesize{$\textbf{M}^{x}$}}} &  \multicolumn{3}{c}{} &  \multicolumn{3}{c}{}  \\
    \cline{1-6}
      & &  &  &  & & & &   \\
      & &  &  &  & & & &  \\
     \multicolumn{3}{c|}{} & \multicolumn{3}{c|}{\raisebox{1\normalbaselineskip}[0pt][0pt]{\footnotesize{$\textbf{M}^{y}$}}} &  \multicolumn{3}{c}{} \\
    \cline{4-9}
     & & \mc{} &  &  & & & &  \\
     & &\mc{}  &  &  & & & &  \\
     \multicolumn{3}{c}{} &  \multicolumn{3}{c|}{}& \multicolumn{3}{c}{\raisebox{1\normalbaselineskip}[0pt][0pt]{\footnotesize{$\textbf{M}^{z}$}}}
  \end{array}
  \right) 
  \renewcommand{\arraystretch}{0.9}
  \left(
  \begin{array}{ c c c | c c c | c c c }
     & &  &  &  & & & &  \\
     & &  &  &  & & & &  \\
    \multicolumn{3}{c|}{\raisebox{1\normalbaselineskip}[0pt][0pt]{\footnotesize{$\textbf{N}^{xx}$}}} &  \multicolumn{3}{c|}{\raisebox{1\normalbaselineskip}[0pt][0pt]{\footnotesize{$\textbf{N}^{xy}$}}} &  \multicolumn{3}{c}{\raisebox{1\normalbaselineskip}[0pt][0pt]{\footnotesize{$\textbf{N}^{xz}$}}}  \\
    \cline{1-9}
      & &  &  &  & & & &   \\
      & &  &  &  & & & &  \\
     \multicolumn{3}{c|}{\raisebox{1\normalbaselineskip}[0pt][0pt]{\footnotesize{$\textbf{N}^{xy}$}}} & \multicolumn{3}{c|}{\raisebox{1\normalbaselineskip}[0pt][0pt]{\footnotesize{$\textbf{N}^{yy}$}}} &  \multicolumn{3}{c}{\raisebox{1\normalbaselineskip}[0pt][0pt]{\footnotesize{$\textbf{N}^{yz}$}}} \\
    \cline{1-9}
     & &  &  &  & & & &  \\
     & &  &  &  & & & &  \\
     \multicolumn{3}{c|}{\raisebox{1\normalbaselineskip}[0pt][0pt]{\footnotesize{$\textbf{N}^{xz}$}}} &  \multicolumn{3}{c|}{\raisebox{1\normalbaselineskip}[0pt][0pt]{\footnotesize{$\textbf{N}^{yz}$}}}& \multicolumn{3}{c}{\raisebox{1\normalbaselineskip}[0pt][0pt]{\footnotesize{$\textbf{N}^{zz}$}}}
  \end{array}
  \right) 
  \right]
    \renewcommand{\arraystretch}{0.9}
    \left(
  \begin{array}{c}
   \\
   \\
  \raisebox{1\normalbaselineskip}[0pt][0pt]{\footnotesize{$\textbf{w}^{x}$}} \\
  \hline
  \\
   \\
  \raisebox{1\normalbaselineskip}[0pt][0pt]{\footnotesize{$\textbf{w}^{y}$}} \\
  \hline
  \\
   \\
  \raisebox{1\normalbaselineskip}[0pt][0pt]{\footnotesize{$\textbf{w}^{z}$}} 
  \end{array}
\right) 
=
j\omega\epsilon_0
  \renewcommand{\arraystretch}{0.9}
 \left(
  \begin{array}{c}
   \\
   \\
  \raisebox{1\normalbaselineskip}[0pt][0pt]{\footnotesize{$\textbf{M}^x\textbf{e}_{\text{inc}}^x$}} \\
  \hline
  \\
   \\
  \raisebox{1\normalbaselineskip}[0pt][0pt]{\footnotesize{$\textbf{M}^y\textbf{e}_{\text{inc}}^y$}} \\
  \hline
  \\
   \\
  \raisebox{1\normalbaselineskip}[0pt][0pt]{\footnotesize{$\textbf{M}^z\textbf{e}_{\text{inc}}^z$}} 
  \end{array}
\right). 
\label{eqn:sys_mat}
\end{equation}
The blocks $\textbf{M}^{x},\ \textbf{M}^y,\ \textbf{M}^z$ are diagonal and each of the blocks $\textbf{N}^{\alpha\beta}$ has block-Toeplitz Toeplitz-block structure on three levels, corresponding to the three physical dimensions of the problem. Note the symmetry in these blocks, i.e., only six of them are unique. Further, each of these blocks is either symmetric or anti-symmetric. This combined with their BTTB structure allows them each to be defined by a single row. Hence the storage cost for the $\textbf{N}$ matrix is $\mathcal{O}(6n)$.

Further note that if the matrix $\textbf{M}$ has a constant diagonal, i.e., the structure is homogeneous, then the matrix $\textbf{I}-\textbf{M}\textbf{N}$ inherits the BTTB structure of $\textbf{N}$. This is the particular case in which circulant preconditioners prove most effective.

\section{A motivating example - uniform waveguide}
\label{sec:straight}
The most elementary and important component in silicon photonics is the straight waveguide, used to channel light signals. From deformations and combinations of waveguides one can create the majority of photonics structures, e.g., Bragg gratings, directional couplers, ring resonators. Therefore, to provide motivation for the preconditioning work to come, and before providing the details of circulant preconditioning, we begin by examining the convergence of GMRES for the unpreconditioned-VIE when simulating guided modes in straight waveguides. In particular, we consider waveguides of various lengths and permittivity contrasts. The most common materials are silicon (Si) and its derivatives: silicon nitride (Si$_3$N$_4$) and silicon dioxide (SiO$_2$), where SiO$_2$ is usually used as the cladding around the waveguide cores. 
The different relative permittivities of typical silicon photonics core and cladding combinations are given in Table~\ref{tab:ref_ind} for reference.
\begin{table}[ht!]
\centering
	\begin{tabular}{c| c}
	Material & Relative permittivity $\epsilon_r$ \\
	\hline
	Silicon (Si)  & 12.1 \\
	Silicon Nitride (Si$_3$N$_4$) & 3.99 \\
	Silicon Dioxide (SiO$_2$) & 2.085 \\
	Si in SiO$_2$ & $12.1/2.085\approx 5.80$ \\
	Si$_3$N$_4$ in SiO$_2$ & $3.99/2.085\approx 1.91 $
	\end{tabular}
	\caption{Relative permittivities of different core and cladding combinations at an operating wavelength of 1550nm. Materials on their own suggest they are cladded by air with $\epsilon_r=1$.}
	\label{tab:ref_ind}
\end{table}

\begin{figure}[hb!]
\centering
\tdplotsetmaincoords{60}{325}
\begin{tikzpicture}[
		scale=0.5,
		tdplot_main_coords,
		axis/.style={->,black,thick},
		cube/.style={very thin,opacity=0.2,fill=red},
		cube hidden/.style={very thin,opacity=.3,gray}]
	\draw[axis] (0,0,0) -- (12.5,0,0) node[anchor=west]{$x$};
	\draw[axis] (0,0,0) -- (0,4,0) node[anchor=east]{$y$};
	
	\draw[cube] (0,0,0) -- (9,0,0) -- (9,0,1) -- (0,0,1) -- cycle;

	\draw[cube] (0,0,0) -- (0,0,1) -- (0,2,1) -- (0,2,0) -- cycle;

	\draw[cube] (0,0,1) -- (0,2,1) -- (9,2,1) -- (9,0,1) -- cycle;
	
	\draw[cube hidden] (9,0,0) -- (9,2,0);
	\draw[cube hidden] (9,2,0) -- (0,2,0);
	\draw[cube hidden] (9,2,0) -- (9,2,1);
	
	\draw[axis] (0,0,0) -- (0,0,3) node[anchor=west]{$z$};

	\foreach \x in {0,0.01,...,2}
		{
 			 \draw[draw=red,opacity=0.1+ 0.22*\x^2]  (9+\x,0,1) -- (9+\x,2,1);
		}

	\foreach \x in {0,0.01,...,2}
		{
 				 \draw[draw=red,opacity=0.1+ 0.22*\x^2]  (9+\x,0,0) -- (9+\x,0,1);
		}
  
	\draw[cube hidden] (11,0,0) -- (11,2,0);
	\draw[cube hidden] (11,2,0) -- (9,2,0);
	\draw[cube hidden] (11,2,0) -- (11,2,1);
	
	\draw (10.1,-0.25,-0.1) node[below]{Absorber};

	\draw (5,2.1,0.6) node[below]{Core};
	\draw (5,-1,0) node[below]{Cladding};
	
\end{tikzpicture}
\caption{Problem setup for a $y$-polarized Gaussian beam directed along the $x$-axis within a strip waveguide. In our experiments, we consider different core and cladding materials.}
\label{fig:strip}
\end{figure}
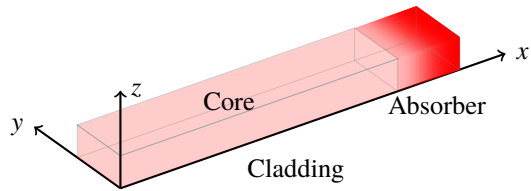

In order to simulate a guided mode in a straight waveguide with VIE, we establish a cuboidal geometry as shown in Fig.~\ref{fig:strip}. The waveguide is excited by a Gaussian beam generated by a dipole placed on the central axis of the waveguide at the left end and directed in the positive $x$-direction. An absorbing region is appended to the right end in order to prevent reflections from the waveguide truncation. We note that, for this simple waveguide structure, only one absorber is required since all waves are propagating in the positive $x$-direction (see \cite{tambova2018adiabatic} for details on domain truncation in VIE).

In the literature one finds many different choices for waveguide cross section. However, for single mode waveguides, their width ($y$-dimension) is slightly larger than one interior wavelength ($\lambda_{\text{int}}$), and their height ($z$-dimension) is slightly larger than half an interior wavelength. For example, for Si-core waveguides, a popular choice is 500nm$\times$220nm~\cite{chrostowski2015silicon}, whereas for Si$_3$N$_4$-core waveguides, a popular choice is 800nm$\times$360nm~\cite{dabos2017tm}. Since we intend to vary refractive index, and hence $\lambda_{\text{int}}$, in this first experiment, we make the following reasonable choice for the waveguide cross-section:
\[
	(\Delta y, \Delta z) = (1.12\lambda_{\text{int}},0.56\lambda_{\text{int}}).
\]
A free-space operating wavelength of 1550nm is considered throughout. 

\begin{figure}[ht!]
		\centering
        \subfigure[Iteration count growing with relative permittivity.]{\includegraphics[width=0.48\textwidth]{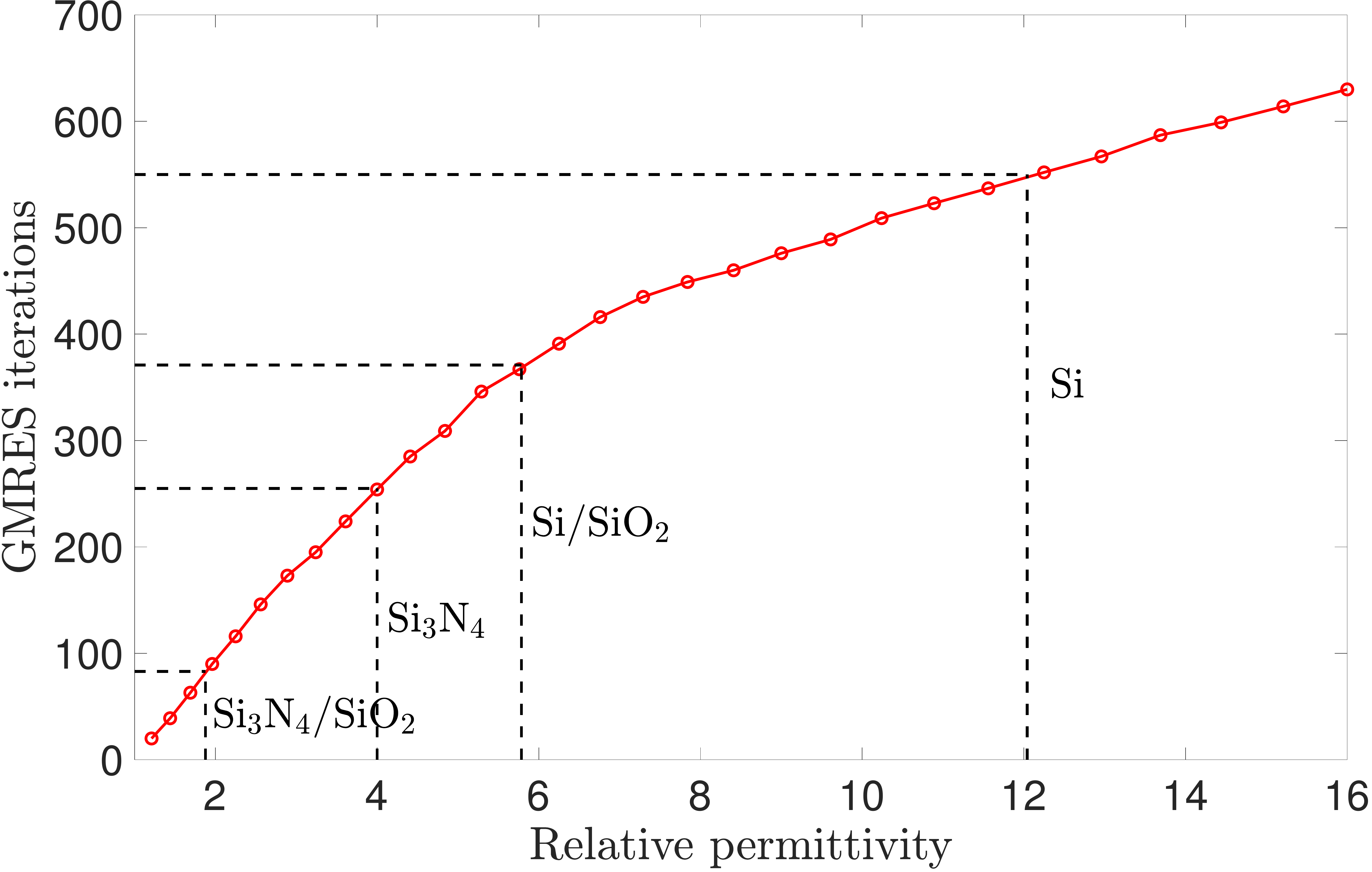}\label{fig:its_v_perm}}
         \subfigure[Iteration count against waveguide length for different core/cladding combinations.]{\includegraphics[width=0.48\textwidth]{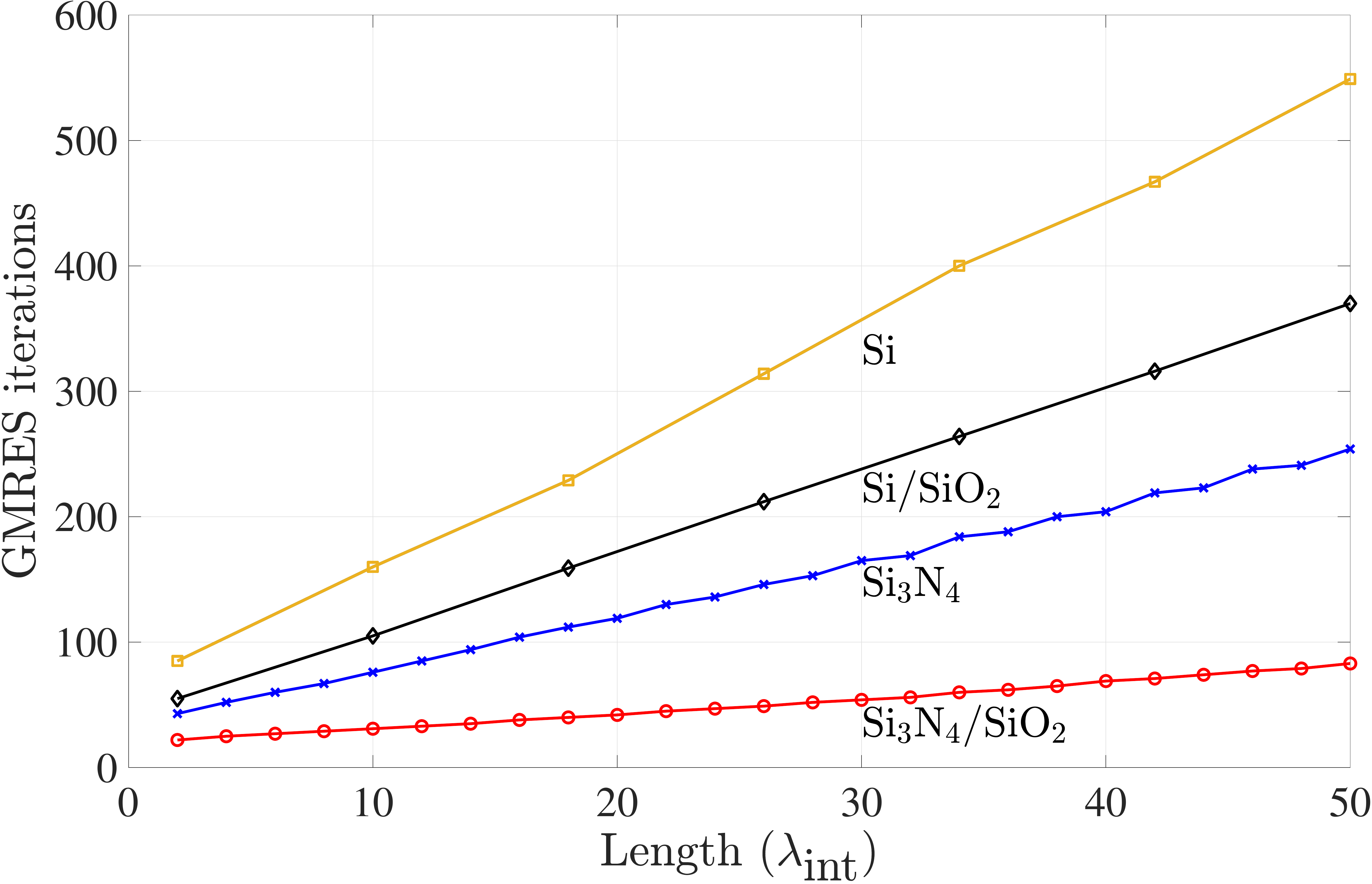}\label{fig:straight_WG_diff_er}}
         \caption{GMRES iteration counts for straight waveguides of different materials and lengths. }
\end{figure}

We begin by fixing the length of the waveguide at 50$\lambda_{\text{int}}$, plus a 5$\lambda_{\text{int}}$ absorber, and increasing the relative permittivity $\epsilon_r$ from 1.2 to 16. For the discretization, a resolution of approximately 20 voxels per $\lambda_{\text{int}}$ is used. The iteration counts required for GMRES convergence to within a tolerance of $10^{-4}$ are shown in Fig.~\ref{fig:its_v_perm}. Observe that the number of iterations required for convergence increases slower than linearly with increasing permittivity, but the increase is significant nonetheless. For example, the iteration count is 83 for silicon nitride with silicon dioxide cladding, compared to 550 for silicon with air cladding.

Next, we perform simulations at a fixed $\epsilon_r$ and  increase $\lambda_{\text{int}}$. The iteration counts for four different core/cladding combinations are shown in Fig.~\ref{fig:straight_WG_diff_er}. We observe that the iteration count grows linearly with waveguide length $\lambda_{\text{int}}$. Compounding this behavior with that observed for increasing permittivity in Fig.~\ref{fig:its_v_perm}, we remark that iteration counts in the hundreds are already observed for waveguides of length $50\lambda_{\text{int}}$ (as seen in Fig.~\ref{fig:straight_WG_diff_er}). 

In the majority of silicon photonics applications, setups with Si cores and SiO$_2$ cladding are used and geometry lengths of up to thousands of wavelengths are considered. As we have demonstrated in this section, such high-constrast high-frequency problems are practically intractable with unpreconditioned-VIE. Therefore, we must devise an appropriate preconditioning strategy. The remainder of this paper discusses an adaptation of an established circulant preconditioning strategy for Toeplitz systems and its application for the first time to electromagnetic problems.  We see that in our silicon photonics setting, circulant preconditioners for VIE can be extremely effective. Furthermore, we present some novel modifications to circulant preconditioners to accelerate them as well as to apply them to the near-Toeplitz systems arising from inhomogeneous structures.

\section{Circulant preconditioning}
\label{sec:prec}
The circulant preconditioners employed here are based on those proposed in \cite{chan1994circulant} for BTTB matrices. 
They are an extension of the optimal point-circulant preconditioners of \cite{chan1988optimal} to the multi-level Toeplitz case. We repeat here the salient features of multi-level circulant preconditioners and refer the reader to \cite{chan1994circulant} for in-depth details.

A Toeplitz matrix $\text{T}_n = \{t_{ij}\}_{i,j=0}^{n-1}$ is Toeplitz if $t_{ij}=t_{i-j}$, i.e., the diagonals are constant. Circulant matrices $\text{C}_n=[c_{ij}]_{i,j=0}^{n-1}$ are also Toeplitz but with the additional property that every row of the matrix is a right cyclic shift of the row above, i.e, $c_{ij}=c_{(i-j)\ \text{mod}\ n}$.  Written out, these matrices have the respective forms
\begin{equation}
\text{T}_n = 
\left( \begin{array}{cccc}
t_0    & t_{-1}  & \ldots &   t_{-(n-1)} \\
t_1    & t_0    &     \ddots       &     \vdots \\
\vdots&   \ddots       &   \ddots    & t_{-1}  \\
t_{n-1}  &  \ldots    & t_1 & t_0   
\end{array} \right)
\end{equation}
and
\begin{equation}
\text{C}_n = 
\left( \begin{array}{cccc}
c_0    & c_{n-1}   & \ldots & c_{1} \\
c_1    & c_0 &   \ddots & \vdots\\
\vdots &   \ddots     &   \ddots  &  c_{n-1}   \\
c_{n-1}  &  \ldots & c_1     & c_0
\end{array} \right).
\end{equation}
We note that circulant matrices have the desirable property that they are diagonalized by the discrete Fourier matrix $\text{F}_n$, such that $\text{C}_n = \text{F}_n^{-1}\Lambda_n\text{F}_n$, where $\Lambda_n=\text{diag}(\text{F}_n\textbf{c})$ is the diagonal matrix consisting of the eigenvalues of $\text{C}_n$. Therefore, $\text{C}_n$ is inverted via the FFT in $\mathcal{O}(n\log n)$ operations. For a Toeplitz matrix, T.\ Chan \cite{chan1988optimal} proposed the \textit{optimal} point-circulant preconditioner whose entries are given by
\begin{equation}
			c_i = 
            \begin{cases}
            				\frac{n-i}{n}t_i + \frac{i}{n}t_{-(n-i)}, & 0\leq i \leq n-1, \\
                            c_{n+i}, 		& -(n-1)\leq i <0.
              \end{cases}
              \label{eqn:Chan}
\end{equation}
This matrix is optimal in the sense that it is the closest circulant matrix to $\text{T}_n$ in the Frobenius norm. There exist other circulant preconditioners (see, for example, the review \cite{strela1996circulant}) and we anticipate the results presented here to be similar for these others. We choose to employ T.\ Chan's preconditioner since it is explicitly defined by the simple formula (\ref{eqn:Chan}) and has been shown to be effective for many Toeplitz problems.

The system matrix for our problem (\ref{eqn:sys_mat}) has BTTB structure. For such matrices, T.\ Chan's preconditioner was extended in \cite{chan1994circulant}. Let us first consider a Toeplitz-block matrix, i.e., one in which each block is point-Toeplitz. Denote such a matrix $\textbf{T}_B$. Then its circulant-block approximation, $\textbf{C}_B$ is obtained by calculating the circulant approximation to each block via (\ref{eqn:Chan}) . If there are $m$ blocks each of size $n\times n$, these matrices are written as  
\[
\textbf{T}_B = 
\left( \begin{array}{cccc}
\text{T}_{11} & \text{T}_{12} & \ldots & \text{T}_{1m} \\
\text{T}_{21} & \text{T}_{22} & \ldots & \text{T}_{2m} \\
\vdots & \vdots  &            & \vdots   \\
\text{T}_{m1} & \text{T}_{m2} & \ldots & \text{T}_{mm}
\end{array} \right)
\]
and
\[
\textbf{C}_B  = 
\left( \begin{array}{cccc}
\text{C}(\text{T}_{11}) & \text{C}(\text{T}_{12}) & \ldots & \text{C}(\text{T}_{1m}) \\
\text{C}(\text{T}_{21}) & \text{C}(\text{T}_{22}) & \ldots & \text{C}(\text{T}_{2m}) \\
\vdots & \vdots  &            & \vdots   \\
\text{C}(\text{T}_{m1}) & \text{C}(\text{T}_{m2}) & \ldots & \text{C}(\text{T}_{mm})
\end{array} \right),
\]
where $\text{C}(\text{T})$ denotes the Chan circulant approximation to $\text{T}$. Having constructed $\textbf{C}_B$, we then proceed to calculate its inverse via applications of the FFT. Each circulant block of $\textbf{C}_B$ has the representation $\text{C}(\text{T}_{ij}) = \text{F}^{-1}\Lambda_{ij}\text{F}$. Defining $\textbf{F} = \text{I}\otimes\text{F}$, we then have that
\begin{equation}
			\textbf{C}_B = [\text{C}(\text{T}_{ij})]_{i,j=1}^m = [\text{F}^{-1}\Lambda_{ij}\text{F}]_{i,j=1}^m = \textbf{F}^{-1}[\Lambda_{ij}]_{i,j=1}^m\textbf{F}.
\end{equation}
The matrix $[\Lambda_{ij}]_{i,j=1}^m$ is an $mn\times mn$ diagonal-block matrix. As described in \cite{chan1994circulant}, this matrix is easily collapsed to a block diagonal matrix $\textbf{D}$ via a permutation matrix $\textbf{P}$, where
\begin{equation}
				\text{diag}(\text{D}_1,\ldots,\text{D}_n) = \textbf{P}[\Lambda_{ij}]_{i,j=1}^m\textbf{P}^{\text{T}}.
				\label{eqn:Di}
\end{equation}
Therefore, the inverse of $\textbf{C}_B$ is given by
\begin{equation}
				\textbf{C}_B^{-1} = \textbf{F}^{-1}\textbf{P}^{\text{T}}\text{diag}(\text{D}_1^{-1},\ldots,\text{D}_n^{-1})\textbf{P}\textbf{F}.
\end{equation}We term $\textbf{C}_B$ the \textit{1-level circulant preconditioner} since we have used one level of circulant approximation. The cost of the inversion of $\textbf{C}_B$ is dominated by the inversion of the $n$ dense blocks $\text{D}_i$, each of size $m\times m$. Therefore, the cost is $\mathcal{O}(nm^3)$. 

If $m$ is small, as is the case for many photonics problems, $\textbf{C}_B$ can be a cheap preconditioner. If $m$ is not small, one may resort to a second level of circulant approximation, applied this time to the dense blocks $\text{D}_i$. In our BTTB case, the blocks $\text{D}_i$ are themselves Toeplitz-block, thus allowing the above procedure to be repeated for each $\text{D}_i$ leading to a \textit{2-level circulant preconditioner} which we denote by $\textbf{C}_{B_2}$. Supposing that each $\text{D}_i$ is comprised of $q$ Toeplitz blocks of size $p\times p$ ($m=pq$), then the cost of inverting $\textbf{C}_{B_2}$ is $\mathcal{O}(npq^3)$. If $q$ is small, $\textbf{C}_{B_2}$ is a cheap preconditioner. For more detailed costings, the reader is again referred to \cite{chan1994circulant}.

In our voxelized setting, we identify the integers $n,p,q$ as $n=n_x$, $p=n_y$, $q=n_z$ described in Section~\ref{sec:VIE} as the numbers of voxels of the discretized domain in the $x, y, z$ directions, respectively. So we can write the costings for the 1-level and 2-level preconditioners as $\mathcal{O}(n_x(n_y n_z)^3)$ and $\mathcal{O}(n_xn_yn_z^3)$, respectively. As we saw in Section~\ref{sec:straight} for the straight waveguide, $n_x\gg n_y,n_z$, and therefore the 1-level preconditioner is feasible. For problems in which $n_y$ is too large, such as disk resonators, we are required to employ the 2-level preconditioner, albeit at the price of reduced performance.

\section{Numerical results}
\label{sec:results}
\subsection{Uniform waveguide}
\label{sec:results_straight}
We return to the uniform waveguide problem examined in Section~\ref{sec:straight}, however without the appended absorber. The omission of the absorber ensures that our domain is homogeneous and thus the diagonal matrix $\textbf{M}$ in (\ref{eqn:discrete_form}) has constant diagonal, and therefore the operator $\textbf{I}-\textbf{M}\textbf{N}$ inherits the BTTB structure of $\textbf{N}$. In sections~\ref{sec:Bragg} and \ref{sec:composite} we consider the case where $\textbf{M}$ has non-constant diagonal.

\begin{figure}[ht!]
\centering
\includegraphics[width=0.5\textwidth]{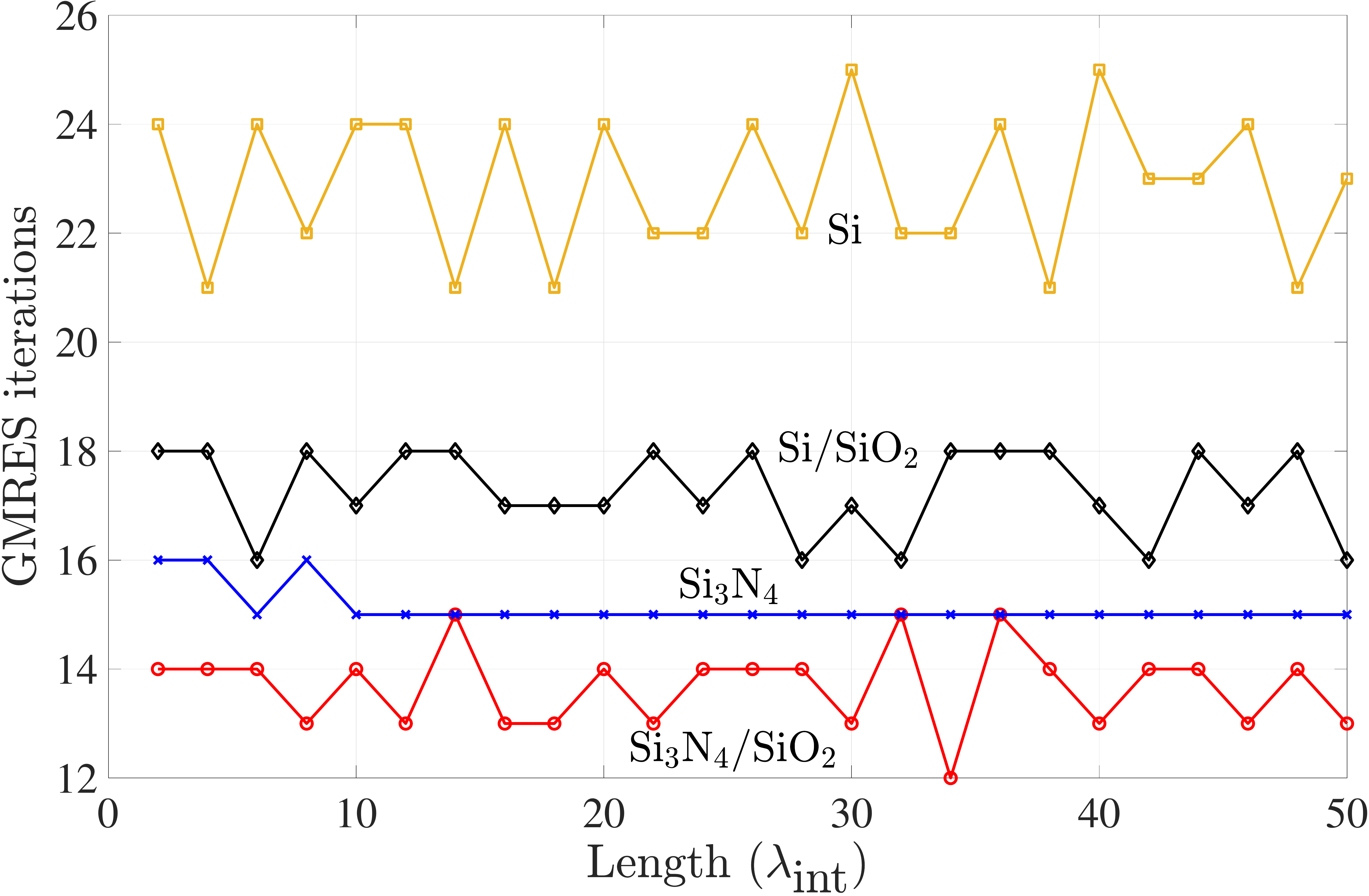}
\caption{GMRES iteration count versus length for straight waveguides of different materials when the 1-level circulant preconditioner is used.}
\label{fig:diff_perms_prec}
\end{figure}

The waveguide has extreme length in only one of the three dimensions, making the 1-level preconditioner described in the previous section cheap to assemble and apply.   Applying the 1-level preconditioner for the waveguide of differing lengths and permittivities allows us to compare directly to Fig.~\ref{fig:straight_WG_diff_er}. The iteration counts are shown in Fig.~\ref{fig:diff_perms_prec}.
Observe that the preconditioner renders the iteration count independent of the waveguide length, although we still observe slightly higher iteration counts for materials with higher refractive index. 

\begin{figure}[ht!]
\centering
		\subfigure[Unpreconditioned.]{\includegraphics[width=0.4\textwidth]{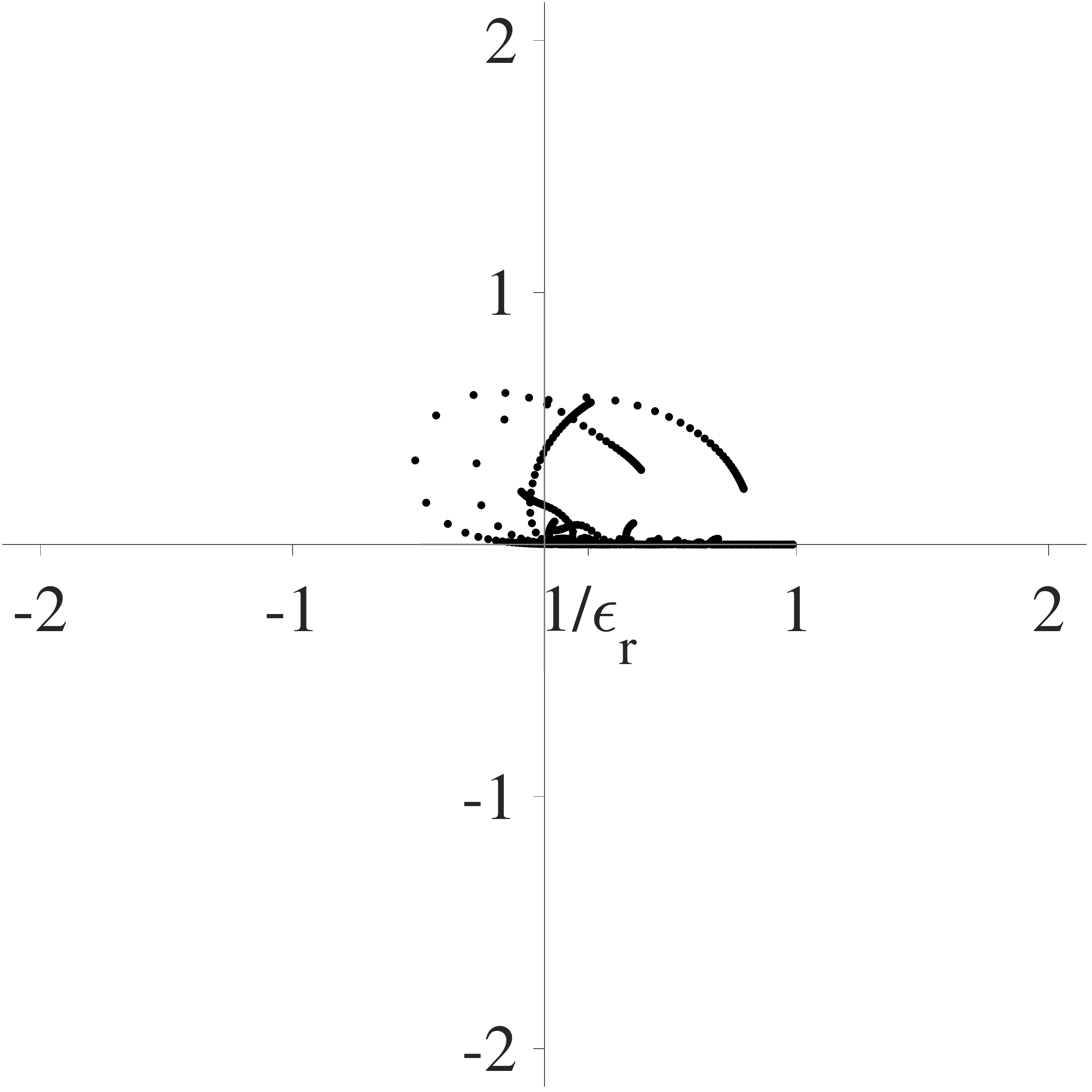}}
		\subfigure[Preconditioned with 1-level circulant.]{\includegraphics[width=0.4\textwidth]{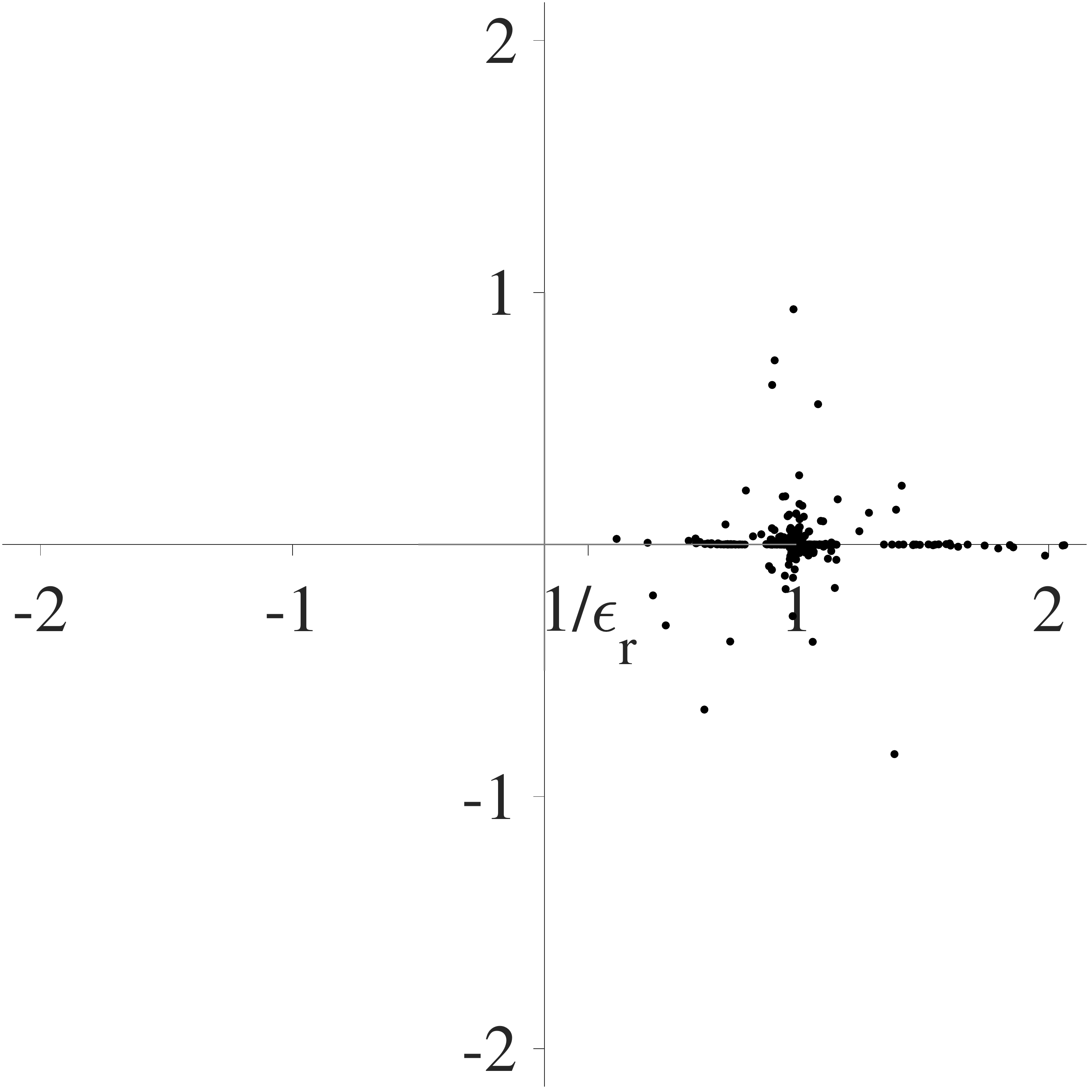}}
\caption{Spectra for Si-core waveguide of length $40\lambda_{\text{int}}$.}
\label{fig:spectra_prec}
\end{figure}

Fig.~\ref{fig:spectra_prec} displays the spectra of the unpreconditioned and preconditioned systems for a silicon waveguide of length $40\lambda_{\text{int}}$. The unprecondioned spectrum is confined to the upper half plane but has many eigenvalues that pass near zero, hence leading to difficulties for iterative solvers. By contrast, the preconditioned spectrum is well separated from zero and clustered around unity, which explains the good performance of iterative solvers for the preconditioned system.

For this simple example, we observe the 1-level preconditioner performs excellently. Next, we consider more difficult scenarios, starting with inhomogeneity in the long $x$-direction. However, before advancing, we make an observation that allows for a dramatic reduction in the cost of the 1-level preconditioner.

\subsubsection{Cost reduction of 1-level circulant preconditioners}
\label{sec:cost_reduc}
In Section~\ref{sec:intro}, we stated that the geometries for which circulant preconditioners are particularly effective are those that are extremely long in one of their three dimensions. However, we saw that the assembly time of the 1-level preconditioner grows as $\mathcal{O}(n_x(n_yn_z)^3)$, and the memory as $\mathcal{O}(n_x(n_yn_z)^2)$. This linear growth in $n_x$ leads to prohibitive cost for extremely long structures (such as Bragg gratings, considered in Section~\ref{sec:Bragg}). However, via a simple observation, we may circumvent this issue without having to resort to 2-level preconditioners, which are cheaper but less effective in reducing iteration counts. 

This observation is that the majority of the $n_x$ blocks, $\text{D}_i$ (see (\ref{eqn:Di})), produced in the 1-level circulant approximation are extremely similar and so may be replaced by a single representative block. The blocks with the more distinct values, which ought to be retained, are those that correspond to $(y,z)$ voxel slices near the two ends of the domain. 

\begin{figure}[h!]
 			\centering
            \subfigure[$\hat{\textbf{v}}$ across blocks $\text{D}_i$ for voxel  $i=1,\ldots,n_x$.]{\includegraphics[width=0.42\textwidth]{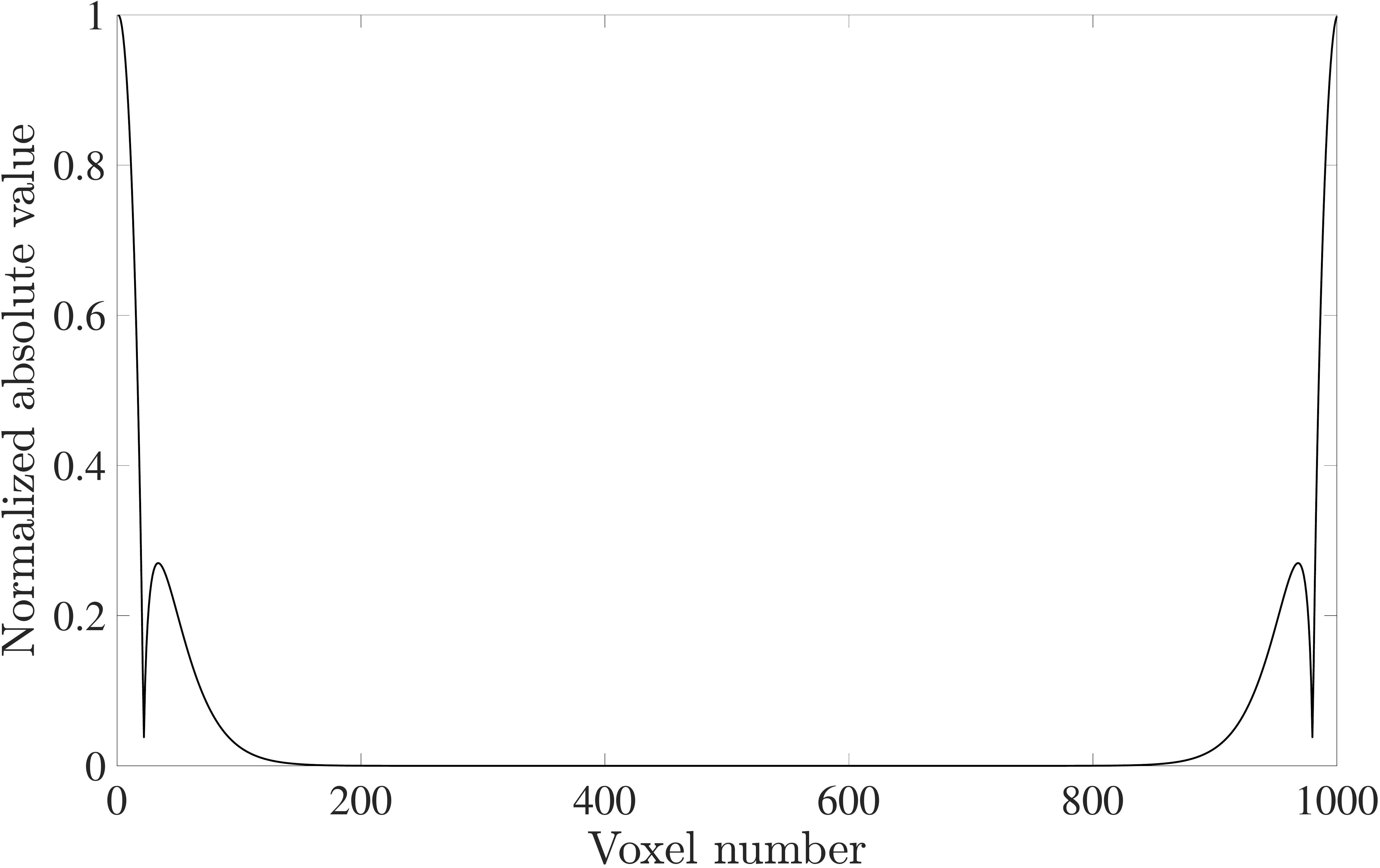}}
            \subfigure[Assembly and solve CPU times.]{\includegraphics[width=0.47\textwidth]{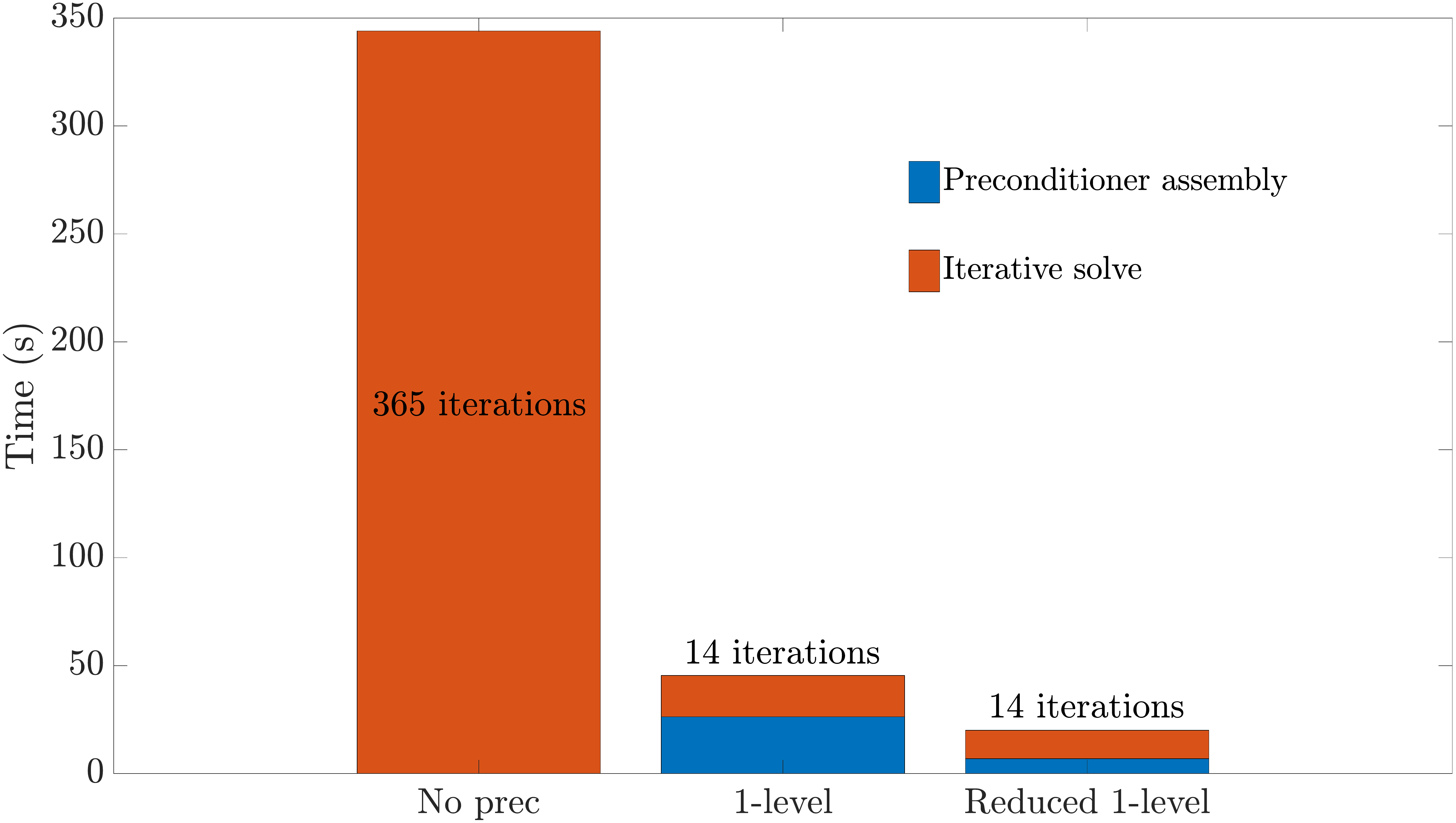}}
            \caption{The reduced 1-level circulant preconditioner for a straight waveguide of length $50\lambda_{\text{int}}$. The values $\hat{\textbf{v}}$ (defined in (\ref{eqn:proxy1})--(\ref{eqn:proxy2})) are used to discard many of the $n_x$ blocks $\left\{\text{D}_i\right\}_{i=1}^{n_x}$: small values lead us to discard the corresponding blocks. This allows for a reduction in setup cost and also solve time for the 1-level preconditioner.}
            \label{fig:abs}
 \end{figure}

We found that a good proxy for determining how many blocks to retain is to look at a representative value from each block $D_i$. Here we choose to consider a vector $\textbf{v}$ of values 
\begin{equation}
    \text{v}_i := \text{D}_i\left(1, \ceil*{\frac{n_y}{2}}\ceil*{\frac{n_z}{2}}\right), \quad \text{for}\ i=1,\ldots,n_x.
    \label{eqn:proxy1}
\end{equation}
More specifically, we consider the relative absolute values of $\textbf{v}$, defined as 
\begin{equation}
    \hat{\textbf{v}} := \left|\frac{\textbf{v}}{||\textbf{v}||_{\infty}}\right|.
    \label{eqn:proxy2}
\end{equation}  
An example $\hat{\textbf{v}}$, for a waveguide of length $50\lambda_{\text{int}}$, is plotted in Fig.~\ref{fig:abs}(a). Observe that near the ends the values of $\hat{\textbf{v}}$ are significant and vary, whereas for a large portion of the domain they are close to zero. By choosing to keep those $\text{D}_i$ where $\hat{\textbf{v}}$ is greater than some tolerance and approximating the remaining $\text{D}_i$ by the central block, $\text{D}_{\ceil{n_x/2}}$ say, we can drastically reduce the computational cost of assembling the preconditioner. In this paper, we use a tolerance of $10^{-3}$, which we found to yield a substantial memory reduction without compromising the iteration count (larger tolerance values yield a greater memory reduction however at the potential expense of an increased iteration count). We remark that this method does not rely on the homogeneity of the structure, as we shall observe in Section~\ref{sec:Bragg} where the method is used for the Bragg grating.

In Fig.~\ref{fig:abs}(b) we show the performance of the cost reduction strategy for a straight silicon waveguide of length $50\lambda_{\text{int}}$. For this example, we keep approximately 12\% of the diagonal blocks in the 1-level preconditioner according to our proxy as defined in (\ref{eqn:proxy1})--(\ref{eqn:proxy2}). This reduces the assembly cost and has the added advantage of accelerating the iterative solve. This acceleration is due to the re-use of the central block in place of the discarded blocks -- since the recycled block stays in the computer's cache, the overhead of data transfer is reduced.
 
To further illustrate the acceleration achieved via this cost reduction strategy, we present in Table~\ref{tab:mem} the memory required to store the 1-level and reduced 1-level preconditioners for waveguides of length $50\lambda_{\text{int}}$ and $500\lambda_{\text{int}}$. The memory required to store the integral operator $\textbf{N}$, defined in (\ref{eqn:op_N}), is provided for reference. Also presented are the times required for matrix-vector products with the respective matrices. We observe that our strategy reduces the memory required to store the preconditioner by approximately a factor of 4 and the average cost of an MVP with the preconditioner is reduced by a factor of approximately 2. This brings down the MVP cost of the 1-level preconditioner to close to the MVP cost of the VIE.
\begin{table}[ht!]
\centering
			\begin{tabular}{c|c|c|c|c}
            & \multicolumn{2}{c|}{50$\lambda_{\text{int}}$} & \multicolumn{2}{c}{500$\lambda_{\text{int}}$} \\
            \cline{2-5}
            & Mem.(GB) & MVP(s) & Mem.(GB) & MVP(s) \\
            \hline\hline
            $\textbf{N}$ operator & 0.236 & 0.296 & 2.36 & 4.01 \\
            1-level   & 12.2 & 0.722 & 121 & 9.25 \\
            Red. 1-level & 2.97 & 0.328 & 29.2 & 3.92 \\
            \end{tabular}
            \caption{Memory requirements in gigabytes and MVP times in seconds for the $\textbf{N}$ operator, the 1-level preconditioner, and the reduced 1-level preconditioner.}
            \label{tab:mem}
 \end{table}

\subsection{Bragg grating -- modifying circulant preconditioning for inhomogeneous structures}
\label{sec:Bragg}
A geometry of practical interest in silicon photonics is that of a Bragg grating. The Bragg grating, as depicted in Fig.~\ref{fig:Bragg}, is a waveguide with a periodically modulating width. Bragg gratings can be up to thousands of wavelengths long, hence can be challenging for numerical methods \cite{chrostowski2015silicon}. However, we observe that their geometry is well-suited for circulant preconditioners since it is close to homogeneous and is long in only one dimension. In order to simulate this device with VIE, we place a Gaussian beam source at the left-hand end of the grating in order to excite the quasi-TE mode. The structure is truncated at either end by adiabatic absorbers~\cite{tambova2018adiabatic}. The geometrical parameters used (shown in Fig.~\ref{fig:Bragg}) here are the following: $d=220\text{nm},\ W=500\text{nm},\ \Delta W = 40\text{nm},\ 
	\Lambda = 320\text{nm},\ N_{\text{per}}=\text{various}$.
\def\nn{15}
\begin{figure}[ht!]
\centering
\begin{tikzpicture}[scale=0.25]

\shade[left color=white,right color=black] (\nn+7,0) rectangle (\nn+9,2);
\shade[left color=black,right color=white] (-4,2) rectangle (-2,0);

\draw (1,2) -- (1,2.25) -- (2,2.25) -- (2,1.75);
\draw (1,0) -- (1,-0.25)-- (2,-0.25)-- (2,0.25);

\foreach \i in {1,3,5,7,9,11,13,\nn}
{
	\draw (\i+1,1.75)--(\i+2,1.75)--(\i+2,2.25)--(\i+3,2.25)--(\i+3,1.75);
	\draw (\i+1,0.25)--(\i+2,0.25)--(\i+2,-0.25)--(\i+3,-0.25)--(\i+3,0.25);
}

\draw (\nn+3,1.75) -- (\nn+4,1.75) -- (\nn+4,2) ;
\draw (\nn+3,0.25) -- (\nn+4,0.25) -- (\nn+4,0) ;

\draw (\nn+4,0) -- (\nn+9,0) -- (\nn+9,2) -- (\nn+4,2);
\draw (1,0) -- (-4,0) -- (-4,2) -- (1,2);

\draw [|-|] (5,2.75) -- (7,2.75);
\draw (6,2.9) node[anchor=south]{$\Lambda$};

\draw [|-|] (1,-0.75) -- (\nn+4,-0.75);
\draw (10,-1) node[anchor=north]{$N_{\text{per}}\times\Lambda$};

\draw [->] (12.5,1.0)--(12.5,1.7);
\draw [->] (12.5,3.0)--(12.5,2.3);
\draw (12.5,2.8) node[anchor=west]{$\Delta W$};

\draw [<->] (-0.75,0.1)--(-0.75,1.9);
\draw (-0.75,1) node[anchor=west]{$W$};

\draw[->] (\nn+8,-1)node[anchor=north]{Absorber}--(\nn+8,-0.1);

\end{tikzpicture}

\caption{Top view of the layout for Bragg grating with period $\Lambda$, width $W$, corrugation depth $\Delta W$, and length $N_{\text{per}}\times\Lambda$, where $N_{\text{per}}$ is an integer. The grating is  truncated with adiabatic absorbers to eliminate reflections. This layout is extruded a distance $d$ out of the page.}
\label{fig:Bragg}
\end{figure}
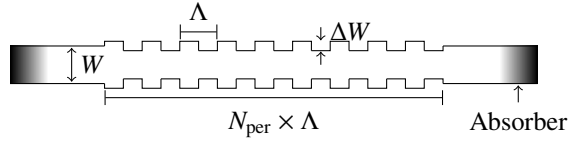

This geometry differs from the straight waveguide in that it is no longer homogeneous in the $x$-direction. Therefore the diagonal matrix $\textbf{M}$ does not have a constant diagonal since the computation domain includes air voxels; furthermore we include the absorbers as shown in Fig.~\ref{fig:Bragg} so at the ends the permittivity is complex. Therefore, the operator $\textbf{I}-\textbf{M}\textbf{N}$ is no longer Toeplitz. In order to re-use the efficient preconditioner from the previous section, we consider a homogenized version of the Bragg grating, with permittivity matrix $\tilde{\textbf{M}}$ that has constant diagonal. Then the corresponding operator $\textbf{I}-\tilde{\textbf{M}}\textbf{N}$ is Toeplitz and we can construct the 1-level circulant approximation $\tilde{\textbf{C}}_B$ to  $\textbf{I}-\tilde{\textbf{M}}\textbf{N}$ as described in Section~\ref{sec:prec}. 

This section investigates, via numerical experiments, an effective way to construct $\tilde{\textbf{M}}\approx\textbf{M}$ such that $\text{diag}(\tilde{\textbf{M}})$ is constant and $\tilde{\textbf{C}}_B$ is an effective preconditioner for $\textbf{I}-\textbf{M}\textbf{N}$. We consider three different ways to construct $\tilde{\textbf{M}}$. First we highlight that, since silicon materials are isotropic, the blocks $\textbf{M}^{\alpha}$, for $\alpha=x,y,z$, are all equal. So we can restrict out discussion to the construction of an individual block $\tilde{\textbf{M}}^{\alpha}$. In order to describe the construction of $\tilde{\textbf{M}}^{\alpha}$, it is convenient to order the $n_xn_yn_z$ entries of $\textbf{M}^{\alpha}$ in an array $\texttt{M}^{\alpha}$ of size $n_x\times n_y\times n_z$ so that $\texttt{M}^{\alpha}(k,l,m)$ contains the value of $(\epsilon_r(\bx)-1)/\epsilon_r(\bx)$ for the voxel in the $k$th position in the $x$-direction, the $l$th position in the $y$-direction, and $m$th position in $z$-direction. The three versions of $\tilde{\texttt{M}}^{\alpha}$ we consider are then given as:  
\begin{enumerate}
	\item $\tilde{\texttt{M}}^{\alpha}(:,:,:) = \text{mode}(\texttt{M}^{\alpha}(:,:,:))$, i.e., the modal average of all the entries in $\texttt{M}^{\alpha}$.
    \item $\tilde{\texttt{M}}^{\alpha}(i,:,:) = \frac{1}{n}\sum_{j=1}^n\texttt{M}^{\alpha}(j,:,:)$ for $i=1,\ldots,n$, i.e., the arithmetic mean along the $x$-direction.
    \item $\tilde{\texttt{M}}^{\alpha}(i,:,:) = \frac{1}{n}\text{Re}\left(\sum_{j=1}^n\texttt{M}^{\alpha}(j,:,:)\right)$ for $i=1,\ldots,n$, i.e., the real part of the arithmetic mean along the $x$-direction.
\end{enumerate}

We perform experiments to assess the performance of these three variations as both the length of the Bragg grating and the modulation depth $\Delta W$ are increased. 

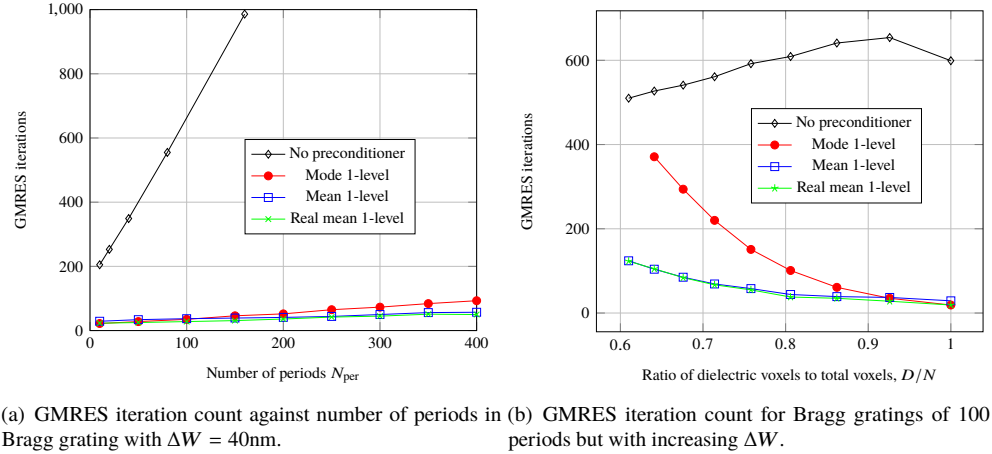
\begin{figure}[h!]
\centering
\subfigure[GMRES iteration count against number of periods in Bragg grating with $\Delta W=40$nm.]{
\begin{tikzpicture}[scale=0.75]
\begin{axis}[
	xlabel={Number of periods $N_{\text{per}}$},
	ylabel={GMRES iterations},
	ymax={1e3},
	ymin={0},
	xmin={0},
	xmax={400},
	grid=both,
	font=\footnotesize,
	legend style={at={(0.4,0.3)},anchor=south west}
]
	
\addplot[color=black,mark=diamond] coordinates {
	(10,205) (20,253) (40,349) (80,555) (160,986)
}; 
		
\addplot[color=red,mark=*] coordinates {
	(10,22) (50,28) (100,35) (150,46) (200,52) (250,65) (300,73) (350,84)
	(400,93)
}; 	

\addplot[color=blue,mark=square] coordinates {
	(10,29) (50,34) (100,37) (150,39) (200,41) (250,44) (300,50)
    (350,56) (400,57)
 };
    
\addplot[color=green,mark=x] coordinates {
	(10,23) (50,25) (100,28) (150,31) (200,36) (250,42) (300,45)
    (350,51) (400,50)
	
}; 
\legend{No preconditioner, Mode 1-level, Mean 1-level, Real mean 1-level}
\end{axis}
\end{tikzpicture}
\label{fig:Bragg_prec}
}
\subfigure[GMRES iteration count for Bragg gratings of 100 periods but with increasing $\Delta W$.]{
\begin{tikzpicture}[scale=0.75]
\begin{axis}[
	xlabel={Ratio of dielectric voxels to total voxels, $D/N$},
	ylabel={GMRES iterations},
	grid=both,
	font=\footnotesize,
	legend style={at={(0.4,0.4)},anchor=south west}
]
	
\addplot[color=black,mark=diamond] coordinates {
	(1,599) (0.926,654) (0.862,641) (0.806,609) (0.758,592) (0.714,561) (0.676,541) (0.641,527) (0.610,510)
}; 
	
\addplot[color=red,mark=*] coordinates {
	(1,19) ( 0.926,35) (0.862,61) (0.806,101) (0.758,151) (0.714,220) (0.676,294) (0.641,371)
}; 

\addplot[color=blue,mark=square] coordinates {
	(1,29) ( 0.926,37) (0.862,39) (0.806,44) (0.758,58) (0.714,69) (0.676,85) (0.641,104) (0.610,124)
}; 

\addplot[color=green,mark=star] coordinates {
	(1,19) ( 0.926,28) (0.862,35) (0.806,38) (0.758,55) (0.714,67) (0.676,84) (0.641,105) (0.610,123)
}; 
			
\legend{No preconditioner, Mode 1-level, Mean 1-level, Real mean 1-level}
\end{axis}
\end{tikzpicture}
}
\caption{Performance of preconditioner for Bragg gratings of various lengths and modulation depths $\Delta W$.}
\label{fig:Bragg_plots}
\end{figure}

First, we examine increasing the length, which is indicated by the number of periods, $N_{\text{per}}$. For a fixed modulation depth $\Delta W=40$nm, $N_{\text{per}}$ is increased from 10 up to 320, which corresponds to a length of 14.2$\mu$m (31.8$\lambda_{\text{int}}$) to 113$\mu$m (254$\lambda_{\text{int}}$). The performance of the three preconditioner types is shown in Fig.~\ref{fig:Bragg_plots}(a). All three preconditioner variations perform excellently, with the iteration count growing very mildly with increasing $N_{\text{per}}$. The best performance is achieved by the third variation, i.e., using the real of the arithmetic mean of the permittivities. For this preconditioner, the iteration count is smaller than 50 for all lengths considered. 

Next, we examine the effect of increasing the size of $\Delta W$ on the performance of the three preconditioner variants. This is an interesting question since, as $\Delta W$ increases, $\tilde{\texttt{M}}$ becomes an increasingly less accurate approximation to $\mathbf{M}$, and so  we might expect the preconditioner to prove ineffective for very large $\Delta W$. The preconditioner performance results are shown in Fig.~\ref{fig:Bragg_plots}(b) for a grating with $N_{\text{per}}$ and $\Delta W$ ranging from 0nm to 280nm. In the figure, the ratio of dielectric (silicon) voxels to total voxels, $D/N$, is displayed rather than $\Delta W$ since it is a quantity more easily generalized to other types of structure. We see that the two arithmetic mean preconditioners still yield significant reductions in iteration counts, even when only 60\% of the volume is dielectric. For most typical Bragg gratings, $\Delta W$ is smaller than 100nm, which translates to $D/N\approx 0.82$ in the figure. Therefore, we can assert that the third variation of the 1-level preconditioner considered is very successful for practical Bragg grating simulations with VIE. 

We remark that this homogenization approach is also well suited to simulations for analyzing manufacturing variations in photonics devices. For example, for assessing the impact of surface wall roughness~\cite{lee2014low}. Typically, these variations are small and lead to the alteration of a small fraction of the voxel permittivities.

\subsection{Composite structures}
\label{sec:composite}
So far we have seen that the 1-level circulant preconditioner performs excellently for straight waveguides and Bragg gratings. However, many silicon photonics problems are composed of multiple components, some of which are not straight and not small in two of their three dimensions. Here we focus on two examples: a directional coupler and a disk resonator. These two devices are composed of two disjoint pieces for which we propose using a \textit{blocked-circulant} preconditioner. In this setup, the geometry is first subdivided into boxes, as shown in Fig.~\ref{fig:disk_15452} and Fig.~\ref{fig:coupler_layout}, and then the relevant circulant preconditioner (1-level or 2-level) is constructed for each box. These circulant preconditioners then form the diagonal blocks of our blocked-circulant preconditioner.

\subsubsection{Disk resonator}
\begin{figure}[hb!]
\centering
	\includegraphics[width=0.8\textwidth,trim={1cm 2cm 1cm 6cm},clip]{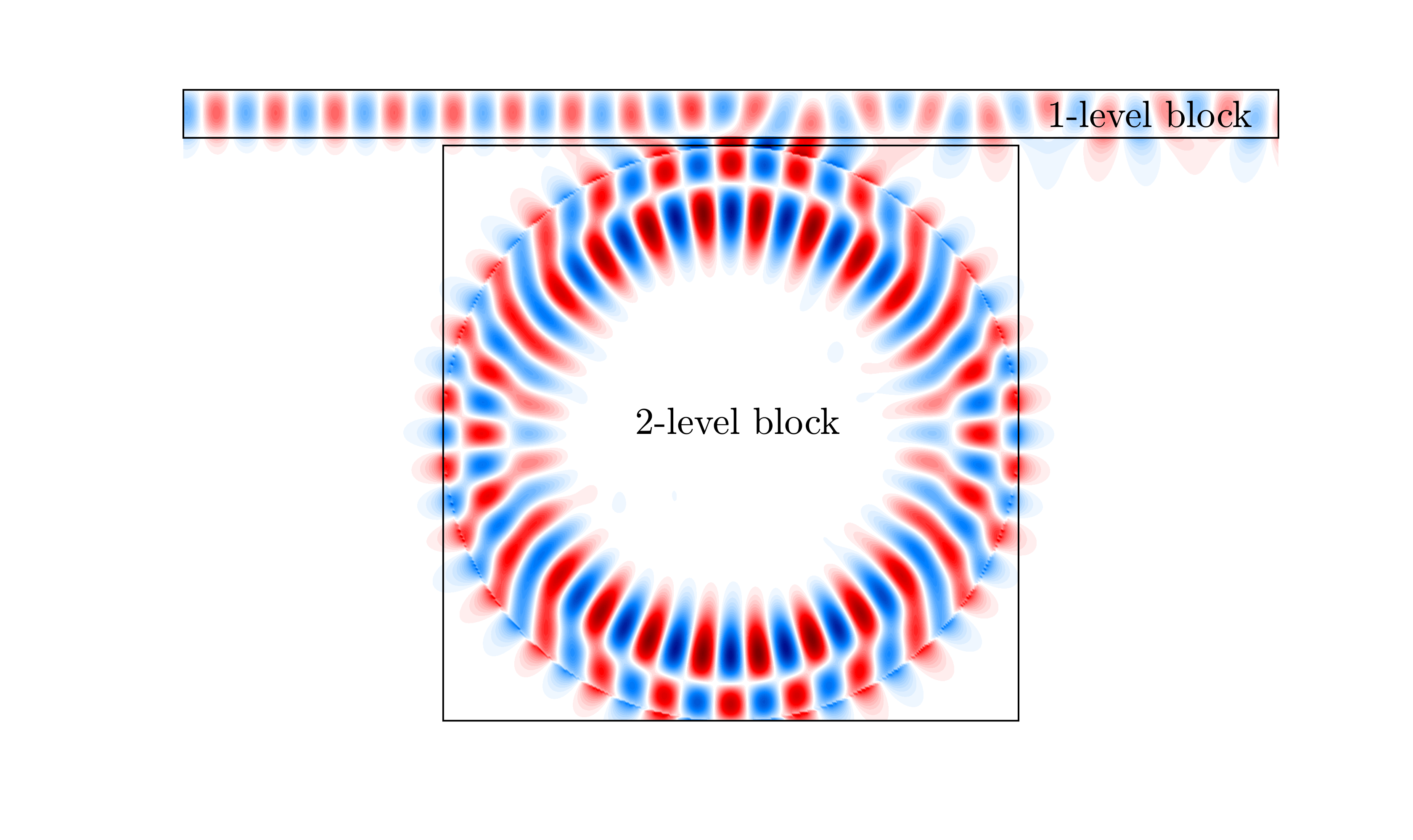}
	\caption{Field plot for a disk resonator of radius 3$\mu$m and for $\lambda_{\text{ext}}=1545$nm. The geometry is partitioned for the construction of the blocked-circulant preconditioner. The 2-level circulant approximation is used for the disk and a 1-level is used for the waveguide.}
	\label{fig:disk_15452}
\end{figure}

Disk resonators, such as that shown in Fig.~\ref{fig:disk_15452}, are important devices in photonics since they behave as spectral filters, which are useful for optical communications applications such as wavelength division multiplexing \cite{bogaerts2012silicon}. 

To construct the blocked-circulant preconditioner, we subdivide the structure into two boxes as shown in Fig.~\ref{fig:disk_15452}. The first box surrounds the bus waveguide for which we construct the reduced 1-level preconditioner as before. The second box encloses the disk. Typical disks have radii up to 5$\mu$m, which corresponds to approximately 11 interior wavelengths. Using voxels of size 20nm (approximately 22 voxels per $\lambda_{\text{int}}$), our discretization will therefore lead to this box having $n_x=242, n_y=242, n_z=11$ voxels in the $x,y,z$ dimensions, respectively (see notation of Section~\ref{sec:VIE}), and therefore a memory requirement of $\mathcal{O}(242(242\cdot 11)^2)$ which is prohibitive on most computers. Therefore, we employ a 2-level circulant approximation for this block, with memory requirement of $\mathcal{O}(242^2 11^2)$, i.e., two orders of magnitude less than the 1-level preconditioner. The value of $\texttt{M}$ used in the construction of the 2-level preconditioner was taken to be the modal average since we found this to yield good performance.

We perform simulations for this disk resonator with disks of radii $1\mu m$-$5\mu m$, with and without the preconditoner. The results are shown in Table~\ref{tab:disk}. We observe that the preconditioner leads to a factor of 22 reduction in iteration count for the largest disk, which equates to a speed-up factor of over 60 times. Although the performance of the preconditioner deteriorates as the disk size is increased, the deterioration is mild and the increase in iteration count is much slower than without a preconditioner.  
\begin{table}[h!]
\centering
	\begin{tabular}{c | c | c | c | c | c }
		\multirow{2}{*}{Radius}		 & \multicolumn{2}{c}{No prec.} & \multicolumn{3}{|c}{Blocked-circulant prec.} \\
	\cline{2-6}
	   ($\mu$m)  & Its. & Solve(s) &  Its. & Build(s) & Solve(s)  \\
	\hline\hline
	1 &  306   & 549    & 40 & 5.41 & 88.6 \\  
	2 &  618   & 3,150  & 83 & 41.6 & 322  \\
	3 &  1,009 & 10,700 & 96 & 141  & 572  \\
	4 &  1,247 & 24,300 & 107 & 320 & 900  \\
	5 &  3,378 & 124,000 & 155 & 689 & 1,940 
	\end{tabular}
	\caption{Simulations for the disk resonator depicted in Fig.~\ref{fig:disk_15452} for disks of increasing radius with 20nm voxels. Shown are the GMRES (tolerance $10^{-4}$) iteration counts with and without preconditioning, the solve times, and the preconditioner build time.}
	\label{tab:disk}
\end{table}

\subsubsection{Directional coupler}
\begin{figure}[h!]
\centering
	\includegraphics[width=0.8\textwidth,trim={15cm 25cm 10cm 25cm},clip]{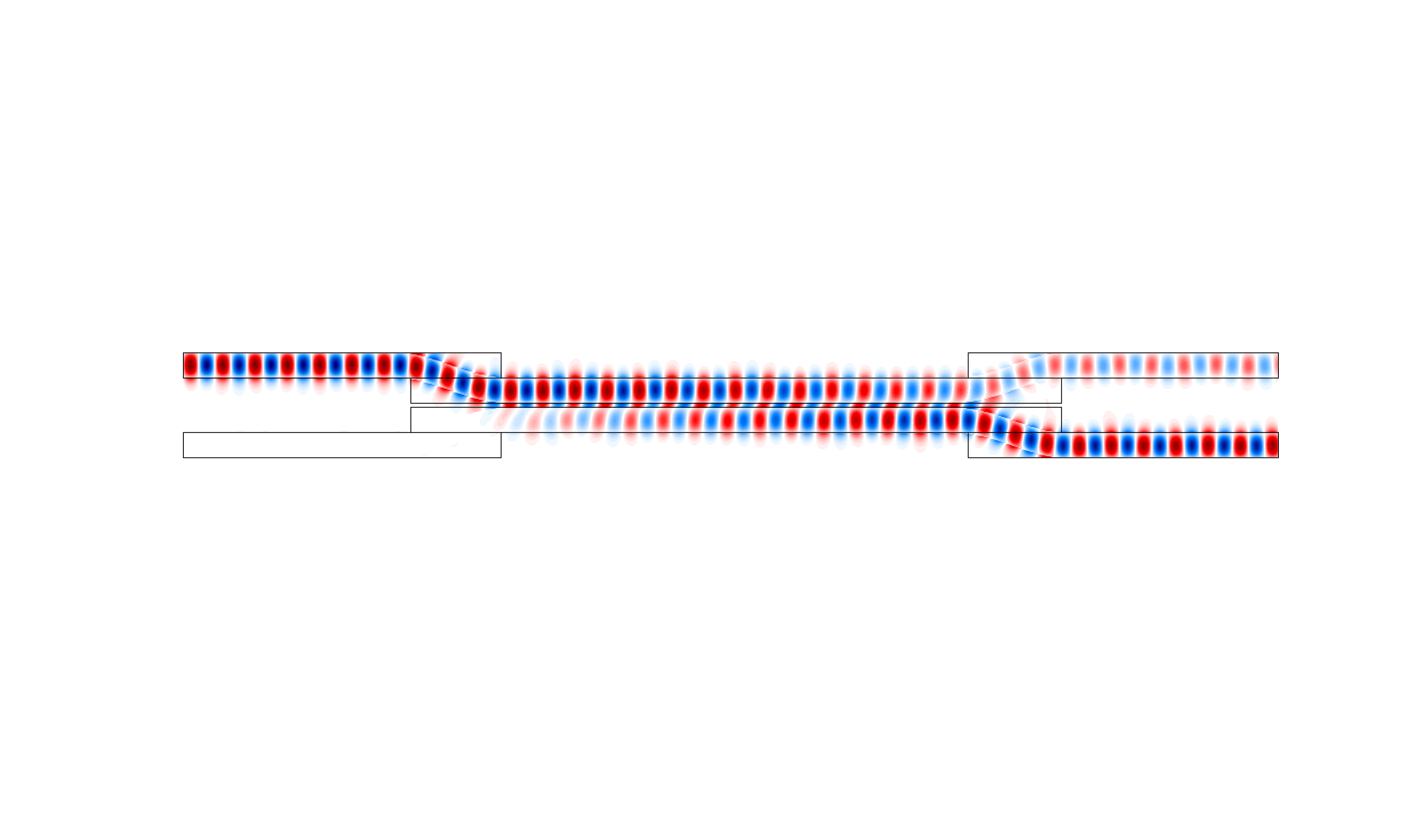}
	\caption{Field plot for a directional coupler. The boxes indicate the partitioning of the geometry for the construction of the blocked-circulant preconditioner.}
	\label{fig:coupler_layout}
\end{figure}

For this example, we consider the blocking strategy shown in Fig.~\ref{fig:coupler_layout}. Since each of the blocks has the width of a waveguide, we may efficiently use 1-level circulant approximations for each of the blocks. 

To test the efficacy of the blocked-circulant preconditioner, we perform simulations, with and without the preconditioner, for directional couplers in which the length $L$ of each of the straight portions is the same, and increases. The results of this experiment are presented in Table~\ref{tab:coupler}. We observe that the iteration count without preconditioner increases linearly in $L$ and is large even for $L=5\lambda_{\text{int}}$. With the preconditioner, the iteration count also increases with $L$, however very slowly: from 77 for $L=5\lambda_{\text{int}}$ to 85 for $L=20\lambda_{\text{int}}$. Compare this to the iteration counts without preconditioner: 321 for $L=5\lambda_{\text{int}}$ to 767 to $L=20\lambda_{\text{int}}$. For the largest problem considered, the preconditioner provides a speed-up factor of more than 10.  
\begin{table*}[h!]
\centering
	\begin{tabular}{c | c | c | c | c | c  | c | c}
		$L$ 	 & \multicolumn{2}{c|}{No prec.} & \multicolumn{3}{c|}{Blocked-circulant prec.} & \multirow{2}{*}{\#voxels} & \#non-air \\
	\cline{2-6}
	   $(\lambda_{\text{int}})$ &   Its.   & Solve(s) & Its. & Build(s) & Solve(s) &  & voxels \\
	\hline\hline
	5 & 321 & 493 & 77 & 21.0 & 175 & 902,055 & 429,550 \\
	10 & 461 & 1,290 & 81 & 28.6 & 270 & 1,288,980 & 613,800 \\
	15 & 607 & 2,680 & 83 & 35.4 & 340 & 1,674,750 & 797,500 \\
	20 & 767 & 5,160 & 85 & 42.3 & 416 & 2,061,675 & 981,750 \\
	\end{tabular}
	\caption{Performance of 1-level blocked-circulant preconditioner for the directional coupler with different waveguide lengths.}
	\label{tab:coupler}
\end{table*}

\section{Conclusion}
\label{sec:conc}
Circulant preconditioners can be extremely effective for the systems of the form $\textbf{I} - \textbf{M}\textbf{T}$ arising in 3D computational electromagnetics. Even though circulant-type preconditioners for multi-level Toeplitz matrices have been proven to not be superlinear~\cite{capizzano2000any}, we have demonstrated here that they can still perform extremely well for the particular geometries arising in silicon photonics applications. The distinguishing characteristic of many of these geometries is that they have extreme length in only one of their three dimensions. 

For the straight waveguide, the 1-level circulant is ideal in that the iteration count is low and independent of the length of the structure. Proving mathematically why this preconditioner is so effective in this simple case would be an interesting avenue for future work. The cost of applying the full 1-level preconditioner is approximately twice that of an MVP with the discretized integral operator $\textbf{N}$. However, the memory required to store the 1-level preconditioner is close to fifty times that required to store $\textbf{N}$. We showed that this storage factor can be reduced to ten times by exploiting the fact that much of the circulant matrix is extremely similar. In addition to reducing the memory footprint of the preconditioner, we saw that this also leads to a speed-up in the preconditioner's application, bringing the per iteration application cost down to being comparable with an MVP with the integral operator. We anticipate that more significant compression of the preconditioner may be possible with further investigation. 

We next considered applying the 1-level circulant preconditioner to the Bragg grating, which has a periodic modulation along its long dimension. Using a permittivity averaging approach we were able to apply the 1-level circulant preconditioner effectively to Bragg gratings that were composed of almost 40\% air voxels. We saw that the iteration count is again low but this time grows very slowly with the length of the grating.

For more complex photonics structures composed of multiple pieces and/or waveguide bends, we proposed a simple blocked-circulant preconditioner based on first partitioning the geometry into boxes and then constructing circulant preconditioners for each box. We saw that for disk resonators of different sizes and for directional couplers of various lengths, this proved an effective preconditioning strategy leading to speed-up factor of over 60 for the most challenging problems. Investigating various geometry partitioning strategies in the future would be useful in order to develop an algorithm to assemble robust and effective blocked-circulant preconditioners for general geometries.

\section*{Funding}
This work was supported by a grant from Skoltech as part of the Skoltech- MIT Next Generation Program, and the Design for Manufacturability (DFM) Methods, PDK Extensions, and Tools for Photonic Systems project sponsored by AIM Photonics.


\bibliography{MIT_papers}

\end{document}